\documentclass[journal]{IEEEtran}
\usepackage{graphicx}
\usepackage{amsmath}
\usepackage{color}
\usepackage[subnum]{cases}
\usepackage{filecontents}
\usepackage{algorithm}
\usepackage{algpseudocode}
\usepackage{float}
\usepackage{setspace}
\usepackage{geometry}
\allowdisplaybreaks
\usepackage{cite}
\usepackage{marvosym}
\usepackage{tabularx}
\usepackage{multicol}
\usepackage{multirow}
\usepackage{textcomp}
\usepackage{epstopdf} %converting to PDF
\usepackage{caption}

\newtheorem{definition}{Definition}

\newtheorem{proposition}{Proposition}

\newtheorem{remark}{Remark}

\newcommand{\NNcg}{\mathcal{N}^{\rm cg}_i}
\newcommand{\NNwg}{\mathcal{N}^{\rm wg}_i}

\newcommand{\qs}{q^{\rm s}_{itw}}
\newcommand{\qDis}{q^{\rm dis}_{itw}}
\newcommand{\qCh}{q^{\rm ch}_{itw}}
\newcommand{\qDisTwo}{q^{\rm dis}_{ikw}}
\newcommand{\qChTwo}{q^{\rm ch}_{ikw}}
\newcommand{\etaDis}{\eta^{\rm dis}_{i}}
\newcommand{\etaCh}{\eta^{\rm ch}_{i}}

\newcommand{\qcg}{q^{\rm cg}_{nitw}}
\newcommand{\qwg}{q^{\rm wg}_{mitw}}
\newcommand{\qcgg}{q^{\rm cg}_{ni(t-1)w}}
\newcommand{\qT}{q^{\rm tr}_{ijtw}}
\newcommand{\qTT}{q^{\rm tr}_{jitw}}

\newcommand{\PRP}[1]{\ensuremath{\Psi_w}}
\renewcommand{\vec}[1]{\ensuremath{\boldsymbol{#1}}}
\begin{document}

\title{Impact of Optimal Storage Allocation on Price Volatility in Electricity Markets}

\author{Amin~Masoumzadeh,
	%        Tansu~Alpcan,~\IEEEmembership{Fellow,~OSA,}\\         
	%        and~Ehsan~Nekouei,~\IEEEmembership{Life~Fellow,~IEEE}
	Ehsan~Nekouei,~\IEEEmembership{Member,~IEEE,}
	Tansu~Alpcan,~\IEEEmembership{Senior Member,~IEEE,}\\
	and~Deb~Chattopadhyay
\thanks{A. Masoumzadeh, E. Nekouei, and T. Alpcan are with the Department
	of Electrical \& Electronic Engineering, University of Melbourne, Australia. D. Chattopadhyay is with the Power System Planning Group of World bank in Washington DC.  E-mails: amin.masoumzadeh@unimelb.edu.au, ehsan.nekouei@unimelb.edu.au, tansu.alpcan@unimelb.edu.au, dchattopadhyay@worldbank.org}
}

\maketitle

%%%%%%%%%%%%%%%%%%%%%%%%%%%%%%%%%%%%%%%%%%%%%%%%%%%%%%%%%%%%%%%%%%%%%%%%%%%%%%%%%55
\begin{abstract}
Recent studies show that the fast growing expansion of wind power generation may lead to extremely high levels of price volatility in wholesale electricity markets. Storage technologies, regardless of their specific forms e.g. pump-storage hydro, large-scale or distributed batteries, are capable of  alleviating the extreme price volatility levels due to their energy usage time shifting, fast-ramping and price arbitrage capabilities. In this paper, we propose a stochastic bi-level optimization model to find the optimal nodal storage capacities required to achieve a certain price volatility level in a highly volatile electricity market. The decision on storage capacities is made in the upper level problem and the operation of strategic/regulated generation, storage and transmission players is modeled at the lower level problem using an extended Cournot-based stochastic game. The South Australia (SA) electricity market, which has recently experienced high levels of price volatility, is considered as the case study for the proposed storage allocation framework. Our numerical results indicate that 80\% price volatility reduction in SA electricity market can be achieved by installing either 340 MWh regulated storage or 420 MWh strategic storage. In other words, regulated storage firms are more efficient in reducing the price volatility than strategic storage firms.

\end{abstract}


\IEEEpeerreviewmaketitle

\vspace*{-.2 cm} 
%%%%%%%%%%%%%%%%%%%%%%%%%%%%%%%%%%%%%%%%%%%%%%%%%%%%%%%%%%%%%%%%%%%%%%%%%%%%%%%%%55
\section{Introduction}

%\IEEEPARstart{E}{lectricity} 

A high level of intermittent wind generation may result in  high price volatility  in electricity markets \cite{ketterer,wozabal,Woo}. In the long term, extreme levels of price volatility can  lead to undesirable consequences such as bankruptcy of retailers \cite{deng} and market suspension. In a highly volatile electricity market, the participants, such as generators, utility companies and large industrial consumers, are exposed to a high level of financial risk as well as costly risk management strategies \cite{volatility}. In some electricity markets, e.g., Australia's National Electricity Market (NEM), which has experienced high levels of price volatility \cite{windTansu}, the market is suspended if the sum of spot prices over a certain period of time is more than cumulative price threshold (CPT). A highly volatile market is subject to frequent CPT breaches due to the low conventional capacity and high level of wind variability.

%The energy-usage shifting and fast-ramping capabilities suggest storage systems for managing the price volatility and facilitating the integration of intermittent renewable resources into the existing electricity networks.  

%In a nodal electricity market which is under a high wind penetration level, the nodal prices become stochastic processes due to the intermittent nature of wind generation. A high level of wind penetration may result in a high price volatility level in the market which is undesirable. \textcolor{blue}{Although price volatility may translate to earnings for peaking  generators, it may induce high levels of uncertainty and risk in generation and network expansion decisions \cite{wozabal}.}
%Wind variability in an energy-only market such as South Australian National Electricity Market  creates high levels of price volatility \cite{windTansu}. 
%The correlation between intermittent power generation and price volatility is also shown in other electricity markets, \emph{e.g.} in Germany \cite{ketterer,wozabal} and in Texas \cite{Woo}. The energy-usage shifting and fast-ramping capabilities of storage systems can  reduce the price volatility and facilitate the integration of renewable resources in the existing electricity networks. 
%Although the current technology cost of storage systems, especially batteries, are relatively high, the  governors' support besides the eventual future decline of the costs can lead to future electrical networks restructured with large capacities of storage systems.     

The current paper proposes a stochastic optimization framework for finding the required nodal storage capacities in electricity markets with high levels of wind penetration such that the price volatility in the market is kept below a certain level. 
The contributions of this paper are summarized as follows: 

\begin{enumerate}
\item  A bi-level optimization model is proposed to find the optimal nodal storage capacities required for avoiding the extreme price volatility levels in a nodal electricity market.

\item  In the upper level problem, the total storage capacities are minimized subject to a price volatility target constraint in each node and at each time. 

\item  In the lower level problem, the non-cooperative interaction between generation, transmission and storage players in the market is modeled  as a stochastic Cournot-based game with an exponential inverse demand function. Note that the equilibrium prices at the lower level problem are functions of the storage capacities. The operation of storage devices at the lower level problem is modeled without introducing binary variables.

\item  The existence of Nash equilibrium  under the exponential inverse demand function is established for the lower level problem.

%Thirdly, we provide the proof that there is no need to put binary variables to force either of charge or discharge of a storage firm zero at the same time according to rationality assumption, which restrains complicating the model with integer variables. Lastly, an exponential function is replaced for the linear inverse demand function in the Cournot-based lower level problem in order to capture the true effect of storage charge/discharge on market price.      
\end{enumerate}

Under the proposed framework, the size of storage devices at two nodes of South Australia (SA) and Victoria (VIC) in NEM is determined such that the market price volatility is kept below a desired level at all times. The desired level of price volatility  can be determined based on various  criteria such as net revenue earned by the market players, occurrence frequency of undesirable prices, number of CPT breaches, etc \cite{AEMC}. 

The proposed storage allocation framework allows the policy makers and market/system operators to compute the required nodal storage capacities for managing the price volatility level in electricity markets.  Although the current cost of storage systems is relatively high, the support from governments (in the form of subsidies) and the eventual decline of the technology cost can lead to  large scale integration of storage systems in electricity markets.

The rest of the paper is organized as follows. The existing related literature is discussed in Section \ref{sec: related}. The system model and  the proposed bi-level optimization problem are formulated in Section \ref{Sec: SM}. The equilibrium analysis of the lower level problem  and the solution method are presented in Section \ref{sec: sol}. The simulation results are presented in Section \ref{sec: sim}. The conclusion remarks are discussed in Section \ref{sec: con}.

\vspace*{-.2 cm} 
\section{Related Works} \label{sec: related}

The problem of optimal storage operation or storage allocation  for facilitating the integration of intermittent renewable energy generators in electricity networks has been studied in \cite{NanLi,VenkatKrishnan,AsmaeBerrada,WeiQi,MahdiSedghi,LeZheng,JunXiao,YuZheng}, with total cost minimization objective functions, and in \cite{RahulWalawalkar,HamedMohsenian2,HosseinAkhavan,HamedMohsenian1,EhsanNasrolahpour}, with profit maximization goals. However, the price volatility management problem using optimal storage allocation has not been investigated in the literature.

%The literature on the topic of storage effect on price volatility is growing substantially. Optimum operation of storage technologies to alleviate the price volatility due to power generation from intermittent renewables is studied in models of cost minimization (optimum allocation and operation models) \cite{NanLi,VenkatKrishnan,MahdiSedghi,AsmaeBerrada,WeiQi,LeZheng,JunXiao} and profit maximization (game theoretical market models) \cite{RahulWalawalkar,HosseinAkhavan,HamedMohsenian2,YuZheng,HamedMohsenian1,EhsanNasrolahpour}.   

The operation of a storage system is optimized, by minimizing the total operation costs in the network, to facilitate the integration of intermittent renewable resources in power systems in \cite{NanLi}. 
 Minimum (operational/installation) cost storage allocation problem for renewable integrated power systems is studied in \cite{VenkatKrishnan,AsmaeBerrada,WeiQi}  under deterministic wind models, and in \cite{MahdiSedghi} under a stochastic wind model.    
The minimum-cost storage allocation problem is studied in a bi-level problem in \cite{LeZheng,JunXiao}, with the upper and lower levels optimizing the allocation and the operation, respectively.  The paper \cite{YuZheng} investigates the optimal sizing, siting, and operation strategies for a storage system to be installed in a distribution company controlled area. 
We note that these works only study the minimum cost storage allocation or operation problems, and the interplay between the storage firms and other participants in the market has not been investigated in these works.

 The paper \cite{RahulWalawalkar} studies the optimal operation of a storage unit, with a given capacity, which aims to maximize its profit in the market from energy arbitrage and provision of regulation and frequency response services. The paper \cite{HamedMohsenian2} computes the optimal supply and demand bids of a storage unit so as to maximize the storage's profit from energy arbitrage in the day-ahead and the next 24 hour-ahead markets. The paper \cite{HosseinAkhavan} investigates the profit maximization problem for a group of independently-operated investor-owned storage units which offer both energy and reserve  in  both day-ahead and  hour-ahead markets. In these works, the storage firm receives the market price as an exogenous input, i.e. the storage is modeled as a price taker firm due to its small capacity.

%{\color{blue}A storage unit aims to maximize its profit in the market, including the day-ahead market and the next 24 hour-ahead markets, by deciding on its energy arbitrage in \cite{HamedMohsenian2}, its energy and regulation services in \cite{RahulWalawalkar}, its energy and reserve trade in \cite{HosseinAkhavan}, and its optimal capacity and operation in \cite{YuZheng}, considering a set of scenarios for wind power availability. Note that considering small price taker storage capacity in the market, the storage firm receives the market price as an exogenous input in these works.}  
  
%Instead of minimizing the costs, a storage unit maximizes its profit considering wind power generation scenarios in the market, including the day-ahead market and the next 24 hour-ahead markets,  in \cite{RahulWalawalkar} from energy arbitrage and regulation services, in \cite{HosseinAkhavan} from energy and reserve trade, in \cite{HamedMohsenian2} from just battery deployment, and in \cite{YuZheng} from optimal storage capacity and operation. 
%Note that considering small price taker storage capacity in the market, the storage firm receives the market price as an exogenous input in these works. 

%Storage and peak demand response have significant strategic values in curbing unwarranted price volatility in a concentrated market that should be considered. 
The operation of a price maker storage device is optimized  using a bi-level stochastic optimization model, with the lower level clearing the market and the upper level maximizing the storage profit by bidding on price and charge/discharge in \cite{HamedMohsenian1}. The storage size in addition to its operation is optimized in the upper level problem in \cite{EhsanNasrolahpour}  when the lower level problem clears the market. Note that the energy and price bids of market participants other than the storage firm are treated exogenously in these models.  

%Given the electricity prices as exogenous inputs in the model, the profit of a storage unit in the day-ahead market and the next 24 hour-ahead markets considering wind power generation scenarios is maximized in \cite{RahulWalawalkar} from energy arbitrage and regulation service, in \cite{HosseinAkhavan} from energy and reserve trade, in \cite{HamedMohsenian2} from just battery deployment, and in \cite{YuZheng} from optimal storage capacity and operation. Moreover, given the both energy and price bids of other market participants as exogenous input in the model, the operation of large storage devices, which are price maker, is optimized in a bi-level stochastic optimization problem, with the upper level clearing the market and the lower level maximizing the storage profit by bidding on price and energy \cite{HamedMohsenian1}. The Storage size in addition to its operation is optimized in the lower level problem, while the upper level problem clears the market given the exogenous energy and price bids of other market participants in \cite{EhsanNasrolahpour}.     

%We are looking for a stochastic market modeling framework in which the market price besides the  generation and consumption decisions of all market participants are calculated concurrently. Papers \cite{BenjaminHobbs,SonjaWogrin} are Cournot-based competition modeling examples and \cite{HobbsMPEC,JalalKazempour} are Bertrand-based competition modeling examples in electricity markets. However, no storage firm is modeled as a participant in those works. 
In \cite{Ventosa,Schill}, the storage firms are modeled as strategic players in Cournot-based electricity markets. However, they do not study storage sizing problem and the effect of intermittent renewables on the market. Therefore, to the best of our knowledge, the problem of finding optimal storage capacity subject to a price volatility management target in electricity markets has not been  addressed before.

\vspace*{-.2 cm} 
%%%%%%%%%%%%%%%%%%%%%%%%%%%%%%%%%%%%%%%%%%%%%%%%%%%%%%%%%%%%%%%%%%%%%%%%%%%%%%%%%55
\section{System Model}\label{Sec: SM}
Consider a nodal electricity market with ${I}$ nodes. Let $\NNcg$ be the set of classical generators, such as coal and gas power plants, located in  node $i$ and $\NNwg$ be the set of wind firms located in  node $i$. The set of neighboring nodes of node $i$ is denoted by $\mathcal{N}_i$. Since the wind availability is a stochastic parameter, a scenario-based model,  with $N_{\rm w}$ different scenarios, is considered to model the wind availability in the electricity network.  % {\bf May not be necessary} The the number of classical generators, \emph{e.g.}, coal, gas and hydro generators, located at node $i$ is denoted by $\Ncg$, and the number of wind generators at node $i$ is denoted by $\Nwg$. 
The nodal prices in our model are determined by solving a Cournot-based game among all market participants, that is, classical generators, wind firms, storage firms and transmission interconnectors which are introduced in detail in the lower level problem. More precisely, the market price in node $i$ at time $t$ under the wind availability scenario $w$ is given by an exponential function:
\begin{align}\label{Eq: Price}
& P_{itw}\left(\vec{q}_{itw}\right)\!=\! \alpha_{it} e^{\!\!-\beta_{it} \! \left(q^{\rm s}_{itw} + \sum\limits_{m \in  \NNwg} q^{\rm wg}_{mitw} +\sum\limits_{n \in \NNcg } q^{\rm cg}_{nitw}+ \sum\limits_{j\in\mathcal{N}_i}q^{\rm tr}_{ijtw} \right)} 
\end{align}
  where  $\alpha_{it}, \beta_{it}$  are positive real values in the inverse demand function, $q^{\rm cg}_{nitw}$ is the generation strategy of the $n$th classical generator located in node $i$ at time $t$ under  scenario $w$, $q^{\rm wg}_{mitw}$ is the generation strategy of the $m$th wind generator located in node $i$ at time $t$ under  scenario $w$, $q^{\rm s}_{itw}$ is the charge/discharge strategy of the storage firm in node $i$ at time $t$ under scenario $w$, $q^{\rm tr}_{ijtw}$ is the strategy of transmission firm located between node $i$ and node $j$ at time $t$ under scenario $w$.  The collection of strategies of all firms located in node $i$ at time $t$ under the scenario $w$ is denoted by $\vec{q}_{itw}$. % Note that the conventional linear inverse demand function is not capable of capturing the true effect of storage on market price. as explained in Appendix \ref{App2}.
  
%	\begin{remark}
%	An alternative approach for finding the market price is using a Bertrand model in which the market price is set by a clearing engine. However, under a Bertrand formulation, our problem will become a three-level problem which is exceedingly hard to analyze.
%	\end{remark}
	
	%This paper investigates and aims to find the minimum required total storage capacity in the market such that the market price volatility stays within desired limit at each time. 
	In this paper, we propose a bi-level optimization approach for finding the minimum required total storage capacity in the market such that the market price volatility stays within a desired limit at each time.

\vspace*{-.2 cm}	
	
\subsection{Upper-level Problem} 
In the  upper-level optimization problem, we determine the nodal storage capacities such that a price volatility constraint is satisfied in each node at each time. In this paper, the variance of market price is considered as a measure of price volatility. The variance of the market price in node $i$ at time $t$, \emph{i.e.} ${\rm Var}\bigl(P_{itw}\bigr)$, can be written as: 
%\begin{subequations} 
	\begin{align}
	&{\rm Var}\left(P_{itw} \right)=\mathsf{E}_w{\left[\left(P_{itw}\left(\vec{q}_{itw}\right)\right)^2\right]}-\left(\mathsf{E}_w\left[{P\left(\vec{q}_{itw}\right)}\right]\right)^2\nonumber \\
	&\hspace{.3 cm}  =   \sum_w \bigg(P_{itw}\left(\vec{q}_{itw}\right)\bigg)^2 \PRP{w} - \bigg(  \sum_w P_{itw}\left(\vec{q}_{itw}\right) \PRP{w}    \bigg)^2  \label{Var}
	\end{align}
%\end{subequations}
where $\Psi_w$ is the probability of scenario $w$.

The notion of variance quantifies the  \emph{effective} variation range of random variables, i.e. a random variable with a small variance has a smaller effective range of variation when compared with a random variable with a large variance.

Given the price volatility relation (\ref{Var}) based on the Nash Equilibrium (NE) strategy collection of all firms $\vec{q}^{\star}_{itw}$, the upper-level optimization problem is given by:

\begin{subequations} \label{ISO}
\begin{align}
&\min_{\left\{Q_i^{\rm s}\right\}_{i}} \sum_{i=1}^I  Q^{\rm s}_i  \nonumber\\
&\text{s.t.}  \nonumber \\
& Q_i^{\rm s} \ge 0 \quad \forall i  \label{ISO1}\\
& {\rm Var}\left(P_{itw}\left(\vec{q}^\star_{itw}\right) \right) \le \sigma_0^2 \quad \forall i,t \label{ISO3} 
\end{align}  
\end{subequations}
where $Q^{\rm s}_i$ is the storage capacity at node $i$, $P_{itw}\left(\vec{q}^\star_{itw}\right)$ is the market price at node $i$ at time $t$ under the wind availability scenario $w$, and $\sigma_0^2$ is the  price volatility target. The price volatility of the market is defined as the maximum variance of market price, \emph{i.e.} $\max_{it} {\rm Var}(P_{itw}(\vec{q}^\star_{itw}))$.
% we aim to achieve. Note that we define the maximum amount of price variance at all regions and all times as the grid-wide price volatility in the paper. 

\vspace*{-.2 cm} 
\subsection{Lower-level Problem}
In the lower-level problem, the nodal market prices and the NE strategies of firms are obtained by solving an extended stochastic Cournot game between wind generators, storage firms, transmission firms, and classical generators. In our formulation, storage and transmission firms can be either regulated or strategic players. 

%, in which some firms are strategic and some are regulated, as distinguished in Definition \ref{def1}. 

\begin{definition} \label{def1}
A strategic firm decides on its strategies over the operation horizon $\{1,...,{N}_{\rm{T}}\}$ such that its aggregate expected profit, over the operation horizon, is maximized. On the other hand, a regulated firm aims to maximize the net market value, i.e. the social welfare. 
\end{definition}

In what follows,  the variable $\mu$ is used to indicate the associated Lagrange variable with its corresponding constraint in the model.

\subsubsection{Wind Generators}
The NE strategy of the $m$th wind generator in node $i$ is obtained by solving the following optimization problem: 
\begin{subequations} \label{Pwind}
	\begin{align}
	&\hspace{0cm}\max_{\left\{q_{mitw}^{\rm wg}\right\}_{tw} \succeq 0} \sum_w \PRP{w} \sum_{t=1}^{{N}_{\rm T}} P_{itw}\left(\vec{q}_{itw}\right)\qwg\nonumber \\
	&{ \rm s.t.} \nonumber \\
	&  \qwg \le {Q}_{mitw}^{ \rm wg} \quad : \quad \mu^{\rm wg,max}_{mitw} \quad \forall t,w \label{cons3}
	\end{align}
\end{subequations}
where $\qwg$ and $Q^{\rm wg}_{mitw}$ are the generation level and the available wind capacity of the $m$th wind generator located in node $i$ at time $t$ under  scenario $w$. Note that the wind availability changes in time in a stochastic manner, and the wind firm's bids depend on the wind availability. As a result, the nodal prices and decisions of the other firms become stochastic in our model. 

%In order to depict the stochasticity in our model, we define random parameter $\phi_w$ multiplied in wind power availability relation as $Q^{\rm wg}_{mitw}= \bar{Q}^{\rm wg}_{mit} \phi_w$, in which $\bar{Q}^{\rm wg}_{mit}$ is the expectation of wind power availability of the $m$th wind generator located in node $i$ at time $t$.

\subsubsection{Storage Firms} Storage firms benefit from price difference at different times to make profit, i.e. they sell the off-peak stored electricity at higher prices at peak times. The NE strategy of  storage firm located in node $i$ is determined by solving the following optimization problem:
\begin{subequations} \label{Pst}
	\begin{align}
	&\max_{\substack{\left\{q_{itw}^{\rm dis},q_{itw}^{\rm ch}\right\}_{tw} \succeq 0\\,\left\{q_{itw}^{\rm s}\right\}_{tw}}} \sum_w \PRP{w} \sum_{t=1}^{{N}_{{\rm T}}} P_{itw}\left(\vec{q}_{itw}\right) \qs -c_{i}^{\rm s} \left(\qDis+ \right.  \nonumber \\ 
	&\hspace{1 cm} \left. \qCh \right) -\gamma^{{\rm s}}_{i} \left(P_{itw}\left(\vec{q}_{itw}\right) \qs+ \frac{P_{itw}\left(\vec{q}_{itw}\right)}{\beta_{it}}\right) \label{storageU} \\	
	&{\rm s.t.} \nonumber \\
	& \qs=\etaDis \qDis - \frac{\qCh}{\etaCh} \quad : \quad \mu^{\rm s}_{itw} \quad \forall t,w \label{cons41}\\
	& \qDis \le \zeta_i^{{\rm dis}} Q^{\rm s}_i \quad : \quad \mu^{\rm dis,max}_{itw} \quad \forall t,w \label{cons42}\\
	& \qCh \le \zeta_i^{{\rm ch}} Q^{\rm s}_i \quad : \quad \mu^{\rm ch,max}_{itw} \quad \forall t,w  \label{cons43}\\
	& 0 \le \sum_{k=1}^{t} \left(\qChTwo-\qDisTwo\right) \Delta \le Q^{\rm s}_i \smallskip : \smallskip \mu^{\rm s,min}_{itw},\mu^{\rm s,max}_{itw} \quad \forall t,w \label{cons44} 
%	& \sum_{t=1}^{{\rm N}_{{\rm T}}} \qCh-\qDis=0 \quad \forall w \label{IFS}
	\end{align}
\end{subequations}
where $\qDis$ and $\qCh$ are the discharge and charge levels of the storage firm in node $i$ at time $t$ under  scenario $w$, respectively, $c_{i}^{\rm s}$ is the unit operation cost,  $\etaCh$ ,$\etaDis$ are the charging and discharging efficiencies, respectively, and $\qs$ is the net supply/demand of the storage firm in node $i$. The parameter $\zeta_i^{{\rm ch}}$ ($\zeta_i^{{\rm dis}}$) is the percentage of  storage capacity $Q^{\rm s}_i$, which can be charged (discharged) during  time period $\Delta$. It is assumed that the storage devices are initially fully discharged. %, used to bound the discharge rate in (\ref{cons42}) and  charge rate in (\ref{cons43}) respectively. 
The energy  level of the storage device in node $i$ at each time is limited by its capacity $Q_i^{\rm s}$. Note that the nodal market prices depend on the storage capacities, \emph{i.e.} $Q_i^{\rm s}$s, through the constraints \eqref{cons42}-\eqref{cons44}. This dependency allows the market operator to meet the volatility constraint using the optimal values of the storage capacities.

The storage firm in node $i$ acts as a strategic firm in the market if $\gamma^{{\rm s}}_{i}$ is equal to zero and  acts as a regulated firm if $\gamma^{{\rm s}}_{i}$ is equal to one. The difference between regulated and strategic players corresponds to the strategic price impacting capability.   Note that the derivative of objective function of the regulated storage firm $i$ is proportional to $P(\cdot)-c^{\rm s}_i$. This intuitively suggests that a regulated storage firm prefers to reduce the market price  to its operation cost while it discharges.

% Higher prices are observed in the market when players act strategically.  

\begin{proposition}
%Considering the unit operation cost of charging/discharging, $c_{i}^{\rm s}$, with any positive amount, we do not need to use binary variables to force either of charging or discharging of a storage device to be zero at any time. In fact, the operation cost from either of charging or discharging prevents calculating positive levels for the both simultaneously in the model.  
At the NE of the lower level game, each storage firm is either in the charge mode or discharge mode, i.e. the charge and discharge levels of each storage firm cannot be simultaneously positive at the NE.
\end{proposition}
{\it{Proof:}} See Appendix \ref{App1}.%{\bf we need a rigorous mathematical argument here.}

\subsubsection{Classical Generators} 
Classical generators include coal, gas, and nuclear power plants. The NE strategy of $n$th classical generator located in node $i$ is determined by solving the following optimization problem:
\begin{subequations} \label{Pcg}
\begin{align}
&\max_{\left\{q_{nitw}^{\rm cg}\right\}_{tw}\succeq 0} \sum_w  \PRP{w} \sum_{t=1}^{{N}_{{\rm T}}} \left(P_{itw}\left(\vec{q}_{itw}\right)-c_{ni}^{\rm cg} \right) \qcg\nonumber  \\
&{\rm s.t.} \nonumber \\
& \qcg \le {Q}_{ni}^{\rm cg}  \quad : \quad \mu^{\rm cg,max}_{nitw} \quad \forall t,w  \label{cons51}\\
& \qcg - \qcgg \le R^{\rm up}_{ni} \quad : \quad \mu^{\rm cg,up}_{nitw} \quad \forall t,w \label{cons52}\\
& \qcgg - \qcg \le R^{\rm dn}_{ni} \quad : \quad \mu^{\rm cg,dn}_{nitw} \quad \forall t,w \label{cons53}
\end{align}
\end{subequations} 
where $ \qcg $ is the generation level of the $n$th classical generator in node $i$ at time $t$ under  scenario $w$,  ${Q}_{ni}^{\rm cg}$ and $c_{ni}$  are the capacity and the short term marginal cost of the $n$th classical generator in node $i$, respectively. The constraints \eqref{cons52} and \eqref{cons53} ensure that the ramping limitations of the $n$th classical generator in node $i$ are always met.      
\subsubsection{Transmission Firms} 
 The NE strategy of the transmission firm between nodes $i$ and $j$ is determined by solving the following optimization problem:
\begin{subequations} \label{Ptr}
	\begin{align}
	&\!\!\!\max_{\left\{\qTT,\qT\right\}_{tw}} \!\!\sum_w \PRP{w}\sum_{t=1}^{{N}_{{\rm T}}} \left(P_{jtw}\left(\vec{q}_{jtw}\right) \qTT
	 \!+\! P_{itw}\left(\vec{q}_{itw}\right)\qT\right)\nonumber\\
	&\hspace{2cm} \left(1-\gamma^{{\rm tr}}_{ij}\right)
	+ \gamma^{{\rm tr}}_{ij}  \left(  \frac{P_{jtw}\left(\vec{q}_{jtw}\right)}{-\beta_{jt}}   +\frac{P_{itw}\left(\vec{q}_{itw}\right)}{-\beta_{it}} \right)     \nonumber \\		
	&{\rm s.t.} \nonumber \\
	& \qT=-\qTT \quad :  \quad \mu^{\rm tr}_{ijtw} \quad \forall t,w\\
	& -Q_{ij}^{\rm tr} \le \qT \le Q_{ij}^{\rm tr} \quad : \quad \mu^{\rm tr,min}_{ijtw},\mu^{\rm tr,max}_{ijtw} \quad \forall t,w \label{cons61}
	\end{align}
\end{subequations} 
where $\qT$ is the electricity exchange level between nodes $i$ and $j$ at time $t$ under  scenario $w$, and $ Q_{ij}^{\rm tr}$ is the capacity of the transmission line between node $i$ and node $j$. The transmission firm between nodes $i$ and $j$ behaves as a strategic player when $\gamma^{{\rm tr}}_{ij}$ is equal to zero and behaves as a regulated player when $\gamma^{{\rm tr}}_{ij}$ is equal to one. Note that the term $P_{jtw}\left(\vec{q}_{jtw}\right) q^{\rm tr}_{jitw}
+ P_{itw}\left(\vec{q}_{itw}\right)q^{\rm tr}_{ijtw}$ in the objective function of the transmission firm is equal to $\left( P_{jtw}\left(\vec{q}_{jtw}\right) - P_{itw}\left(\vec{q}_{itw}\right) \right) q^{\rm tr}_{jitw}$ which implies that the transmission firm between two nodes  makes profit by transmitting electricity from the node with  lower market price to the node with  higher market price.

Transmission lines or interconnectors are usually controlled by the market operator and are regulated to maximize the social welfare in the market. The markets with regulated transmission firms are discussed as  electricity markets with transmission constraints in the literature, e.g., see  \cite{JudithCardel,WilliamHogan,EvaggelosKardakos}. However, some electricity markets allow the transmission lines to act strategically, i.e. to make revenue by trading electricity across the nodes \cite{AEMO}.

%\begin{proposition}
%	Considering the unit operation cost of import/export, $c_{ij}^{\rm tr}$, with any positive amount, we do not need to use binary variables to force either of import or export of a line to be zero at any time. In fact, the operation cost from either of import or export prevents calculating positive levels for the both simultaneously in the model.  
%\end{proposition}
\vspace*{-.2 cm} 

\section{Solution Approach} \label{sec: sol}
In this section, we first provide a game-theoretic analysis of the lower-level problem. Next, the bi-level price volatility management problem is transformed to a single optimization Mathematical Problem with Equilibrium Constraints (MPEC). 
 
\vspace*{-.2 cm}  
 
\subsection{Game-theoretic Analysis of the Lower-level Problem}
To solve the lower-level problem, we need to study the best response functions of  firms participating in the market. Then, any intersection of the best response functions of all firms will be a NE. In this subsection, we first establish the existence of NE for the lower-level problem. Then, we provide the necessary and sufficient conditions which can be used to solve the lower-level problem. 

To transform the bi-level price volatility management problem to a single level problem, we need to ensure that for every vector of storage capacities, \emph{i.e.} $\vec{Q}^{\rm s}=\left[Q_1^{\rm s},\cdots,Q_I^{\rm s}\right]^\top\geq\vec{0}$, the lower-level problem admits  a NE. At the NE strategy of the lower-level problem, no single firm has any incentive to unilaterally deviate its strategy from its NE strategy.  Note that the objective function of each firm is quasi-concave in its strategy and constraint set of each firm is closed and bounded for all $\vec{Q}^{\rm s}=\left[Q_1^{\rm s},\cdots,Q_I^{\rm s}\right]^\top\geq\vec{0}$. Thus, the lower level game admits a NE. This result is formally stated in Proposition \ref{pro2}.

\begin{proposition} \label{pro2}
For any vector of storage capacities, $\vec{Q}^{\rm s}=\left[Q_1^{\rm s},\cdots,Q_I^{\rm s}\right]^\top\geq\vec{0}$, the lower level game admits a Nash Equilibrium.
 %each firm has an objective function continuous in the strategy space and quasi-concave in its strategy, and its strategy space is non-empty, compact and convex. Therefore, according to Theorem 1.2 in \cite{quasiconvex}, the lower-level game admits a Nash Equilibrium.
\end{proposition} 

\begin{IEEEproof}
	Note that the objective function of each firm is continuous and quasi-concave in its
	strategy. Also, the strategy space is non-empty, compact and
	convex. Therefore, according to Theorem 1.2 in \cite{quasiconvex}, the lower level
	game admits a NE.
\end{IEEEproof} 

\subsubsection{Best responses of wind firm $mi$} 
Let $\vec{q}_{-(mi)}$ be the strategies of all firms in the market except the wind generator $m$ located in node $i$.  Then, the best response of  the wind generator $m$ in node $i$ to $\vec{q}_{-(mi)}$ satisfies the necessary and sufficient Karush-Kuhn-Tucker (KKT) conditions ($t \in \{1, ...,  N_{\rm T}\}; w \in \{1, ..., N_{\rm w}\}$):
\begin{subequations} \label{KKTw}
	\begin{align}
& P_{itw}\left(\vec{q}_{itw}\right)+\frac{\partial{P_{itw}\!\left(\vec{q}_{itw}\right)}}{\partial{\qwg}}\qwg  - \frac{\mu^{\rm wg,max}_{mitw}}{\PRP{w}}\le 0 \perp \qwg \ge 0   \label{8LPs}\\
&  \qwg \le {Q}_{mitw}^{ \rm wg} \perp \mu^{\rm wg,max}_{mitw} \ge 0 \label{9s1}
	\end{align}
\end{subequations} 
where the perpendicularity sign, $\perp$, means that at least one of the adjacent inequalities must be satisfied as an equality \cite{ferris}.
\subsubsection{Best responses of storage firm $i$} 
To study the best response of the storage firm in node $i$, let $\vec{q}_{-i}$ denote the collection of  strategies of all firms except the storage firm in node $i$. Then, the best response of the storage firm in node $i$ is obtained by solving the following KKT conditions ($t \in \{1, ...,  N_{\rm T}\}; w \in \{1, ..., N_{\rm w}\}$): 
\begin{subequations} \label{KKTst}
	\begin{align}
&  P_{itw} \left(\vec{q}_{itw}\right)+(1-\gamma^{{\rm s}}_{i} )\frac{\partial{P_{itw} \left(\vec{q}_{itw}\right) }}{\partial{\qs}}\qs   +\!\frac{\mu^{\rm s}_{itw}}{\PRP{w}}=0 \label{9LPs}\\
& \frac{-\!\etaDis \mu^{\rm s}_{itw}\!-\!\mu^{\rm dis,max}_{itw}\!-\!\Delta\!\sum_{k=t}^{{{N}_{{\rm T}}}}\! \left(\!\mu^{\rm s,min}_{ikw}\!-\!\mu^{\rm s,max}_{ikw} \!\right)}{\PRP{w}}-c_{i}^{\rm s} \le 0 \nonumber\\
& \hspace{6cm} \perp \qDis \ge 0 \\
& \frac{\frac{\mu^{\rm s}_{itw}}{\etaCh}\! +\!\mu^{\rm ch,min}_{itw}\!-\!\mu^{\rm ch,max}_{itw}\!+\!\Delta\!\sum_{k=t}^{{{N}_{{\rm T}}}}\! \left(\!\mu^{\rm s,min}_{ikw}\!-\!\mu^{\rm s,max}_{ikw}\! \right)}{\PRP{w}} \!-\!c_{i}^{\rm s} \! \le \! 0 \nonumber \\
& \hspace{6cm} \perp \qCh \ge 0 \\
& \qs=\etaDis \qDis - \frac{\qCh}{\etaCh} \\
& \qDis \le \zeta_i^{{\rm dis}} Q^{\rm s}_i \perp \mu^{\rm dis,max}_{itw} \ge 0\\
& \qCh \le \zeta_i^{{\rm ch}} Q^{\rm s}_i \perp \mu^{\rm ch,max}_{itw} \ge 0\\
& 0 \le \sum_{k=1}^{t} \left(\qChTwo-\qDisTwo\right) \Delta  \perp  \mu^{\rm s,min}_{itw} \ge 0\\
&  \sum_{k=1}^{t} \left(\qChTwo-\qDisTwo\right) \Delta \le Q^{\rm s}_i \perp  \mu^{\rm s,max}_{itw} \ge 0 \label{10s2}
	\end{align}
\end{subequations}

\subsubsection{Best responses of classical generation firm $ni$} 
The best response of the classical generator $n$ in node $i$ to $\vec{q}_{-(ni)}$, i.e. the collection of strategies of all firms except the classical generator $n$ in node $i$, is obtained by solving the following KKT conditions ($t \in \{1, ...,  N_{\rm T}\}; w \in \{1, ..., N_{\rm w}\}$):
\begin{subequations} \label{KKTcg}
	\begin{align}
&   P_{itw}\left(\vec{q}_{itw}\right)\!-\!c_{ni}^{\rm cg}\!+\! \frac{\partial{P_{itw}\left(\vec{q}_{itw}\right)}}{\partial{\qcg}}\!\qcg    \!+\!\frac{\!-\!\mu^{\rm cg,max}_{nitw}\!+\!\mu^{\rm cg,up}_{ni(t+1)w}}{\PRP{w} } \nonumber \\
& +\frac{-\mu^{\rm cg,up}_{nitw}+\mu^{\rm cg,dn}_{nitw}-\mu^{\rm cg,dn}_{ni(t+1)w}}{\PRP{w} }\le 0 \perp \qcg \ge 0 \label{10LPs}\\
&  \qcg \le {Q}_{ni}^{\rm cg} \perp \mu^{\rm cg,max}_{nitw} \\
& \qcg - \qcgg \le R^{\rm up}_{ni} \perp \mu^{\rm cg,up}_{nitw} \ge 0\\
& \qcgg - \qcg \le R^{\rm dn}_{ni} \perp \mu^{\rm cg,dn}_{nitw} \ge 0 \label{11s2}
	\end{align}
\end{subequations}

\subsubsection{Best responses of transmission firm $ij$} 
Finally, the best response of the transmission firm between nodes $i$ and $j$, to $\vec{q}_{-(ij)}$, i.e. the set of all firms' strategies except those of the transmission line between nodes $i$ and $j$, can be obtained  using the KKT conditions ($t \in \{1, ...,  N_{\rm T}\}; w \in \{1, ..., N_{\rm w}\}$):
\begin{subequations} \label{KKTtr}
	\begin{align}	
& P_{itw}\left(\vec{q}_{itw}\right) +\left(1-\gamma^{{\rm tr}}_{ij}\right) \frac{\partial{P_{itw}\left(\vec{q}_{itw}\right)}}{\partial{\qT}}\qT +
\frac{\mu^{\rm tr}_{jitw}+\mu^{\rm tr}_{ijtw}}{\PRP{w}} \nonumber\\
& \hspace{3.5cm} +\frac{\mu^{\rm tr,min}_{ijtw}-\mu^{\rm tr,max}_{ijtw}}{\PRP{w}}=0  \label{11LPs}\\
& \qT=-\qTT\\
& -Q_{ij}^{\rm tr} \le \qT \perp \mu^{\rm tr,min}_{ijtw}\ge 0\\
& \qT \le Q_{ij}^{\rm tr} \perp \mu^{\rm tr,max}_{ijtw}\ge 0   \label{12s1}
	\end{align}
\end{subequations} 

\vspace*{-.2 cm} 
\subsection{The Equivalent Single-level Problem}

Here, the bi-level price volatility management problem is transformed into a single-level MPEC. To this end, note that for every vector of storage capacities the market price can be obtained by solving the firms' KKT conditions. Thus, by imposing the KKT conditions of all firms as constraints in the optimization problem \eqref{ISO}, the price volatility management problem can be written as the following single-level optimization problem:
%\begin{subequations} 
	\begin{align}
	& \min \sum_{i=1}^I  Q^{\rm s}_i  \label{optimization}\\
	&{\rm s.t.} \nonumber \\ 
	& (\ref{ISO1}-\ref{ISO3}),(\ref{8LPs}-\ref{9s1}),(\ref{9LPs}-\ref{10s2}),(\ref{10LPs}-\ref{11s2}),(\ref{11LPs}-\ref{12s1})  \nonumber\\
	& m \in \{1,..., N^{\rm wg}_i\} , n \in \{1, ..., N^{\rm cg}_i\}, i,j \in \{1,..., I\}, \nonumber \\
	& t \in \{1, ...,  N_{\rm T}\}; w \in \{1, ..., N_{\rm w}\}  \nonumber
	\end{align}
%\end{subequations} 
where the optimization variables are the storage capacities, the bidding strategies of all firms and the set of all Lagrange multipliers. Because of the nonlinear complementary constraints, the feasible region is not necessarily convex or even connected. Therefore, increasing the storage capacities stepwise, we solve the lower level problem, which is convex. Once the price volatility constraint is addressed, the optimum solution is found.

	\begin{remark}
%It is possible to convert a MPEC problem to a Mixed-Integer Linear Problem (MILP), e.g. \cite{HamedMohsenian1}, and solve it by a commercial MILP solver such as CPLEX. However, considering the inverse demand function exponential instead of linear, the first order optimality conditions of our problem are not linear. 
It is possible to convert the equivalent single level problem (\ref{optimization}) to a Mixed-Integer Non-Linear Problem (MINLP). However, the large number of integer variables potentially makes the resulting MINLP computationally infeasible.   
	\end{remark}

\vspace*{-.2 cm} 
 
\section{Case Study and Simulation Results} \label{sec: sim}

In this section, we study the impact of  storage installation on price volatility in two nodes of Australia's National Electricity Market (NEM): South Australia (SA) and Victoria (VIC). SA has a high level of wind penetration and VIC has high coal-fueled classical  generation.  Real data for price and demand from the year 2013
%, which are averaged over the year for a day (24 hours), 
is used to calibrate the inverse demand function in the model. Different types of generation firms, such as coal, gas, hydro, wind and biomass, with generation capacity (intermittent and dispatchable) of 3.7 GW and 11.3 GW were active in SA and VIC, respectively, in 2013.
%(total generation capacity in NEM was 45 GW). 
%Recently two existing coal power plants of Playford Steam Turbine (240 MW) and Northern Brown Coal (530 MW) have been retired in SA \cite{AER}. 
The transmission line interconnecting SA and VIC, which is a regulated line, has the capacity of 680 MW but currently is working with just 70\% of its capacity.  
The generation capacities in our numerical results are gathered from Australian Electricity Market Operator's (AEMO's) website (aemo.com.au) and all the prices are shown in Australian dollar.       

Similar to \cite{Morales}, we consider a scenario based analysis wherein three scenarios, i.e. high wind scenario (with probability of 0.2), low wind scenario (with probability 0.2) and base wind scenario (with probability of 0.6), are defined  to capture the wind power availability. %, as shown in Figure \ref{fig_wind}.
 The base wind scenario indicates the available wind generation level for a day (24 hours), in each node,  averaged over a year \cite{AEMOwind}. 
Given that the wind turbines are dispersed over the whole region in each node, we assume that the wind power availability is often around its expected value, i.e. the base wind level. 
%Wind forecast error is assumed normally distributed and the forecast error is $\phi \%$ above and under its base amount \cite{Ding}.      
The wind generation level at high wind and low wind scenarios are assumed to be $\phi\%$ above and below the wind generation level at the base wind scenario, respectively. Various levels of wind availability can be captured by changing the wind power fluctuation parameter $\phi$ \cite{Ding}.

%\begin{figure}[!htb]  
%	%\captionsetup{justification=centering}
%	\centering
%	\includegraphics[scale= 0.65]{wind.pdf} 
%	\caption{Percentage of hourly available wind power in SA with fluctuation of $2\phi \%$ around its base amount.} 			
%	\label{fig_wind}
%\end{figure}
%

In what follows, by price volatility we mean the maximum variance of market price, i.e. $\max_{it}{\rm Var}(P_{itw}(q^\star_{itw}))$. Also, by square root of price volatility we mean the maximum standard deviation of market price, i.e. $\max_{it}\sqrt{{\rm Var}(P_{itw}(q^\star_{itw})}$.

\vspace*{-.2 cm} 

\subsection{One-node model simulations in South Australia}
\begin{figure}[!htb] 
	%\captionsetup{justification=centering}
	\centering
	\includegraphics[scale=.6]{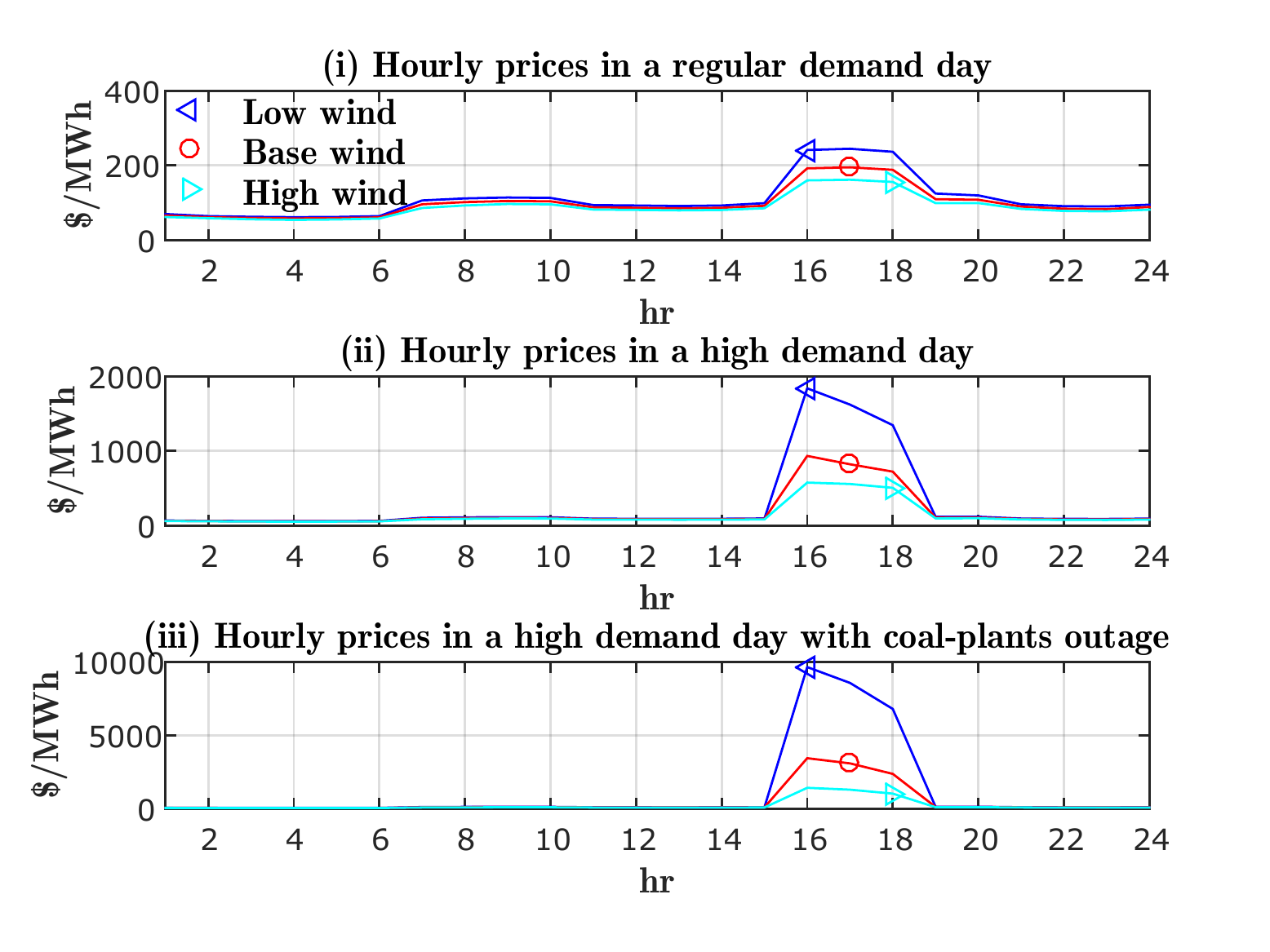} 
	\caption{Hourly wholesale electricity prices in SA with $\phi=50\% $ and no storage.} 			
	\label{fig_price}
\end{figure}

In this subsection, we first study the impacts of peak demand levels and supply capacity shortage on the electricity price in SA with no storage. Next, we study the effect of storage on price volatility in SA. Fig. \ref{fig_price} shows the hourly prices for a day in SA (with no storage) for three different cases: $(i)$ a regular demand day, $(ii)$ a high demand day, $(iii)$ a high demand day with coal plants outage. An additional load of  1000 MW  is considered  in the high demand case during hours 16, 17 and 
18
%, \emph{i.e.} the coefficient $\alpha$ in the inverse demand function is multiplied by $e^{950\beta}$ in those three hours, 
to study the joint effect of wind intermittency  and large demand variations on the price volatility. The additional loads are sometimes demanded in the market due to unexpected high temperatures happening  in the region. 
The coal-plants outage case is motivated by the recent retirement of two coal plants in SA with total capacity of 770 MW \cite{AER}. This allows us to investigate the joint impact of wind indeterminacy and low supply capacity on the price volatility. 

%Note that the square root of price volatility indicates the maximum standard deviation of price, $\max_t (VAR(P_{tw}))^{\frac{1}{2}}$, and has the unit of \$/MWh.

According to Fig. \ref{fig_price}, in a regular demand day, wind power fluctuation with  $\phi=50\%$  does not create much price fluctuation. 
In the regular demand day, the maximum price is equal to 194 \$/MWh in the base wind scenario  whereas it changes to 161 \$/MWh and 244 \$/MWh in the high wind and the low wind scenarios, respectively. The square root of the price volatility  in the regular demand day is equal to 26 \$/MWh.
Based on Fig. \ref{fig_price}, the maximum price in a high demand day in SA changes from 933 \$/MWh in the base wind scenario to 576 \$/MWh and 1837 \$/MWh in the high wind and the low wind scenarios, respectively. The square root of the price volatility in the high demand day is equal to 420 \$/MWh. The extra load at peak times and  the wind power fluctuation create a higher level of price volatility during a high demand day compared with a regular demand day.

The outage of coal plants in SA beside the extra load at peak hours increases the price volatility due to the wind power fluctuation. The maximum price during the high demand day with coal plants outage varies from 3446 \$/MWh to 1436 \$/MWh and 9634 \$/MWh in the high wind and the low wind scenarios, respectively. The square root of the price volatility  during the high demand day with coal plant outage is equal to 2787 \$/MWh. The square root of the price volatility during the high demand day with coal plant outage is almost 107 times more than the regular demand day due to the simultaneous variation in both supply and demand.

%The peak price changes from the base wind case level of 933 \$/MWh between 576 and 1837 \$/MWh in high and low wind cases, with standards deviation of 420 \$/MWh. Moreover, the closure of the coal plants in SA beside the extra load at peak times accelerates the price volatility resulted from the wind power fluctuation. The peak price shifts from the base wind case level of 3446 \$/MWh to 1436 and 9634 \$/MWh, with standard deviation of 2787 \$/MWh. The standard deviation of price volatility becomes almost 107 times more with the both changes in supply and demand. 

%Storage operation reduces the price volatility beside shaving the peak price. The minimum storage capacities required for achieving different price volatility reduction targets in SA are calculated by solving the optimization model (\ref{optimization}) in a high demand day with coal-plants outage case.

\begin{figure}[!htb] 
	%\captionsetup{justification=centering}
	\centering
	\includegraphics[scale=.6]{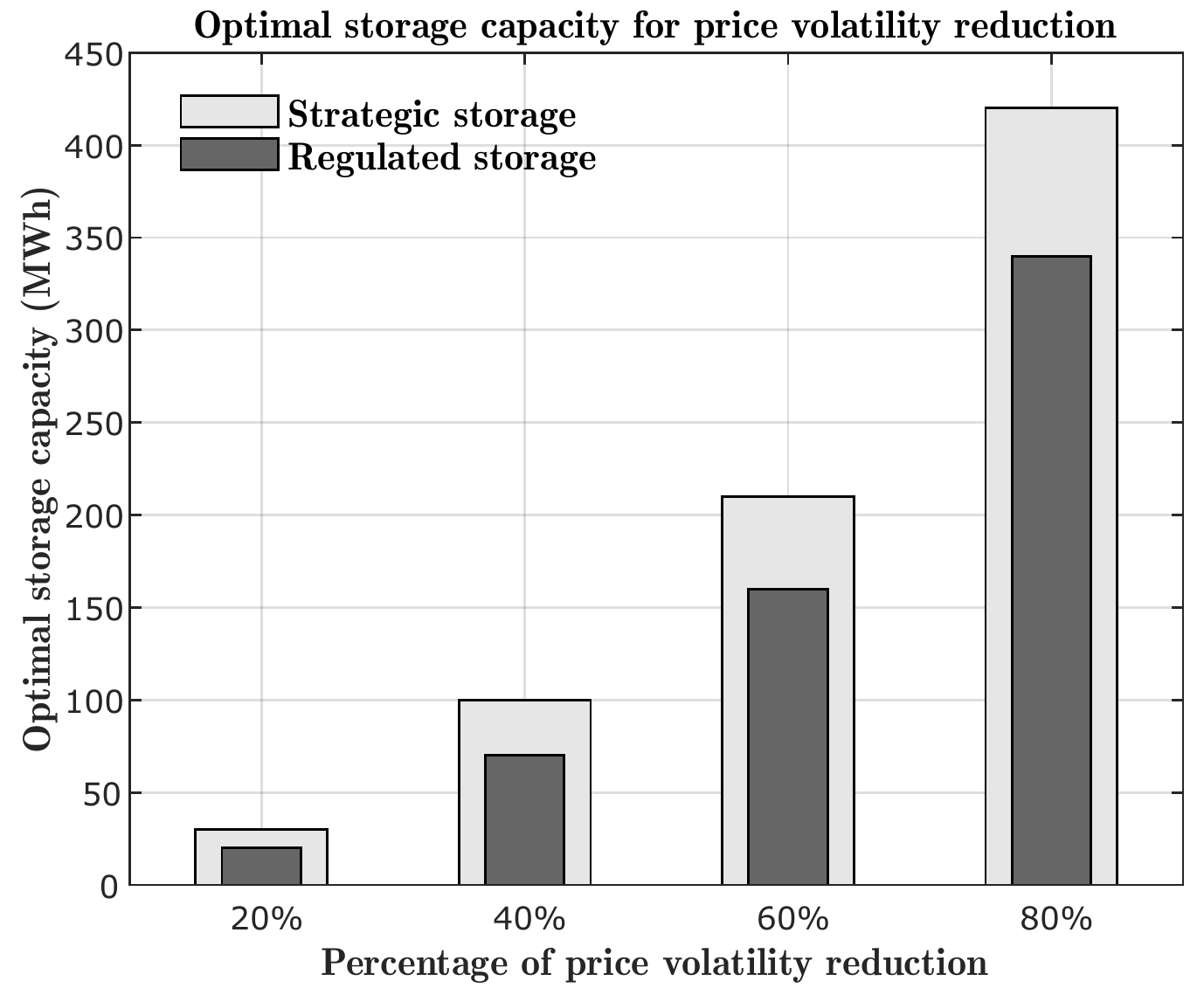} 
	\caption{Optimal strategic and regulated storage capacity for achieving different price volatility levels in SA node with $\phi=50\%$ for a high demand day with coal-plants outage.} 			
	\label{table1}
\end{figure}

Fig. \ref{table1} shows the minimum required (strategic/regulated) storage capacities for achieving various levels of price volatility in SA during a high demand day with coal plants outage. The minimum storage capacities  are calculated by solving the optimization problem (\ref{optimization}) for the high demand day with coal-plants outage case with $\phi=50\%$.
According to Fig. \ref{table1}, a strategic storage firm requires a substantially larger capacity, compared with a regulated storage firm, to achieve a target price volatility level due to the selfish behavior of the storage firms.
In fact, the strategic storage firms may sometimes withhold their available capacities and do not participate in the price volatility reduction as they do not always benefit from reducing the price. 
The price volatility in SA can be reduced by 80\% using either 420 MWh strategic storage or 340 MWh regulated storage.
Note that AEMO has forecasted about 500 $\rm MWh$ battery storage to be installed in SA until 2035 \cite{BattFuture}.

%large enough storage capacity, around the forecasted amount of 500 MWh battery storage to be installed in next 20 years in SA,  in either type of strategic or regulated is capable of abating price volatility at least by 80\%.   

According to our numerical results, storage can displace the peaking generators, with high fuel costs and market power, which results in reducing the price level and the price volatility. A storage capacity of 500 MWh (or $\frac{500}{2}$ MW given the discharge coefficient $\eta^{\rm dis}=\frac{1}{2}$) reduces the square root of the price volatility from 2787 \$/MWh to 919 \$/MWh, almost 30\% reduction, during a high demand day with coal-plant outage in SA.

\begin{figure}[!htb] 
	%\captionsetup{justification=centering}
	\centering
	\includegraphics[scale=0.6]{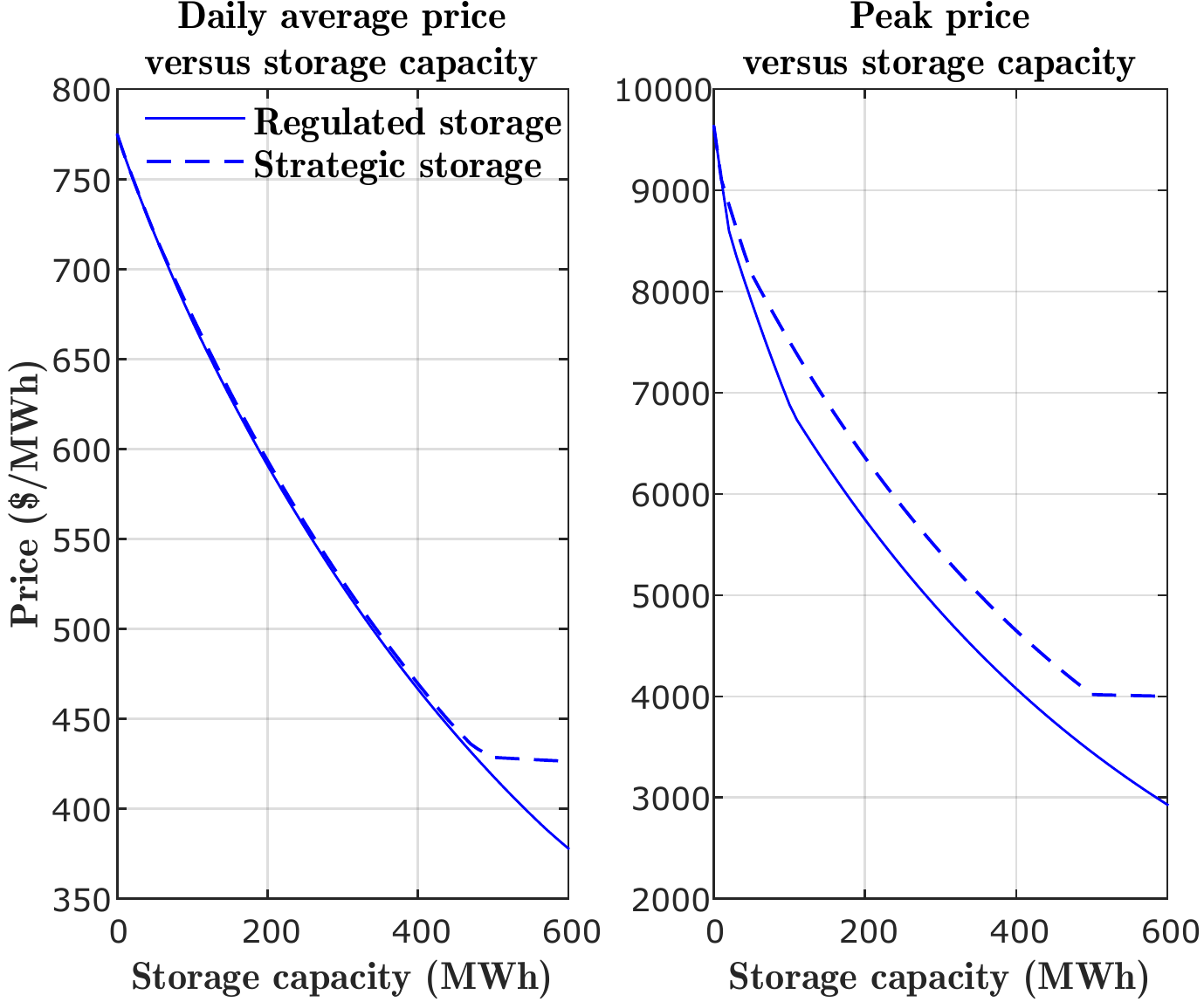} 
	\caption{Daily peak and average prices in SA versus storage capacity with wind power fluctuation parameter $\phi=50\%$ in a high demand day with coal-plant outage.} 			
	\label{fig_AP}
\end{figure}

The behaviour of the peak and the daily average prices for the high demand day with coal plants outage in SA is illustrated in Fig. \ref{fig_AP}. In this figure, the peak price represents the maximum of price over all scenarios during the day, i.e. $\max_{t,w}P_{tw}(\vec{q}^\star_{tw})$ and the daily average price indicates the average of price over time and scenarios, i.e. $\frac{1}{N_T}\sum_{tw}P_{tw}(q^\star_{tw})\Psi_w$.
Sensitivity analysis of the peak and the daily average prices in SA with respect to storage capacity  indicates that high storage capacities lead to relatively low prices in the market. 
At very high prices, demand is almost inelastic and a small amount of excess supply leads to a large amount of price reduction.
According to Fig. \ref{fig_AP}, the rate of price reduction decreases as the storage capacity increases since large storage capacities lead to relatively low peak prices which make the demand more elastic.

Based on Fig. \ref{fig_AP}, the impact of storage  on the daily average and peak prices depends on whether the storage firm is strategic or regulated. It can be observed  that the impacts of strategic and regulated storage firms on the daily peak/average prices are almost similar for small storage capacities,
i.e. when the storage capacity is smaller than 100 MWh (or $\frac{100}{2}$ MW given $\eta^{\rm dis}=\frac{1}{2}$). However, a regulated firm reduces both the peak and the average prices more efficiently compared with a strategic storage firm  as its capacity becomes large. 
%indicates more efficient behavior with larger capacities. 
%In fact, both the strategic and regulated storage firms stop operating when they reach their optimal point of charging/discharging.  
A large strategic storage firm in SA does not use its excess capacity beyond  500 MWh to reduce the market price since it acts as  a strategic profit maximizer, but a regulated storage firm contributes to the price volatility reduction  as long as there is potential for price reduction by its operation.

\begin{figure}[!htb]  
	%\captionsetup{justification=centering}
	\centering
	\includegraphics[scale=0.6]{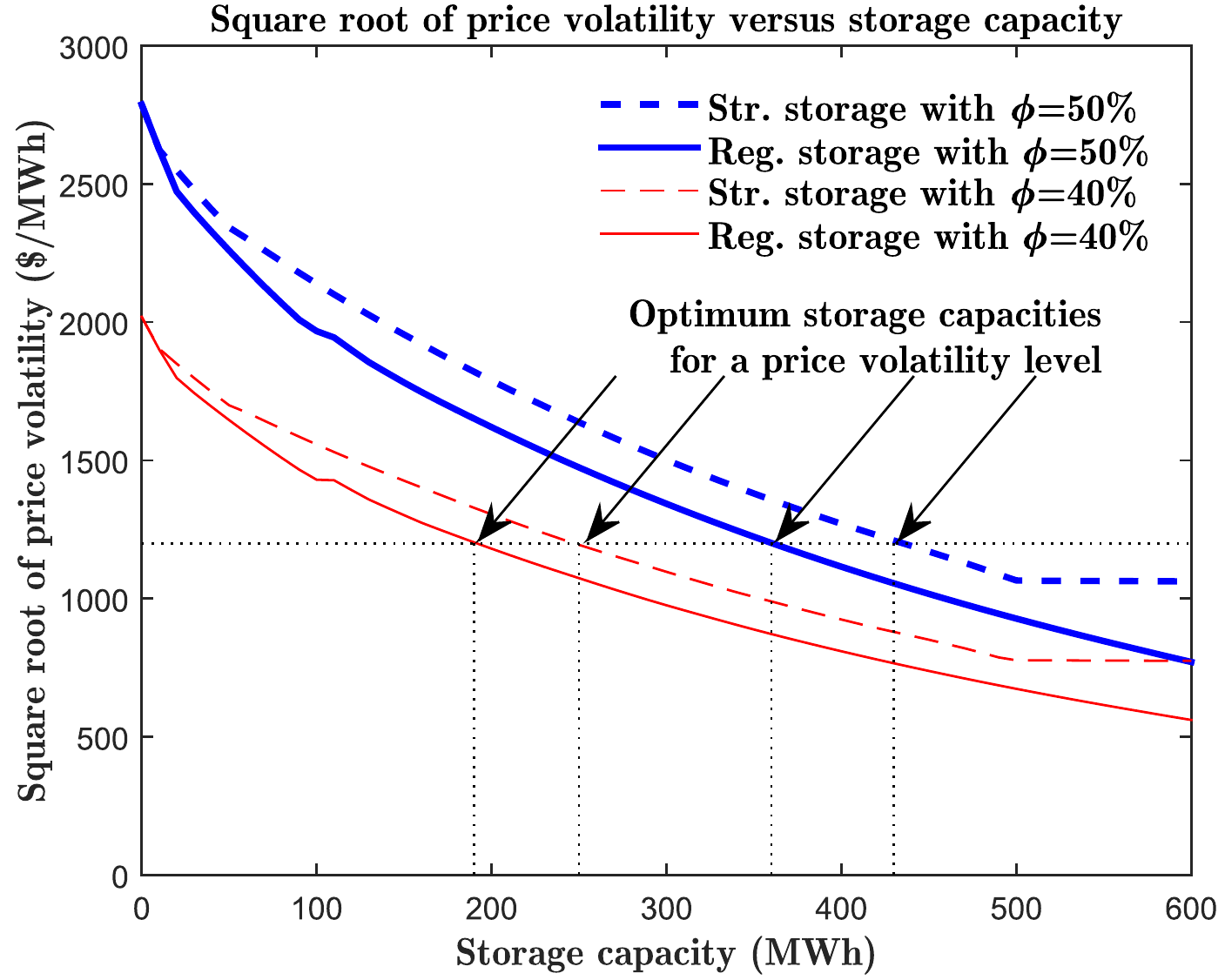} 
	\caption{Square root of price volatility in SA versus storage capacity  with  $\phi$ equal to $40\%$ and $50\%$ during a high demand day with coal-plants outage.} 			
	\label{fig_vol}
\end{figure}

Fig. \ref{fig_vol} depicts the square root of price volatility  in SA during the high demand day with coal plant outage for $\phi=50\%$ and $\phi=40\%$. 
%Mathematically, the maximum standard deviation of price can be written as $\max_t (VAR(P_{tw}))^{\frac{1}{2}}$.
According to this figure, the price volatility in the market decreases by installing either regulated or strategic storage devices.
However, a strategic storage firm stops reducing the price volatility  when its capacity exceeds a threshold value. 
%Moreover, by decreasing the wind power fluctuation parameter $\phi$ from 50\% to 40\%, the square root of price volatility  diminishes almost by 30\%.
 Moreover, the square root of price volatility, in all cases, diminishes almost by 30\% as the wind
 power fluctuation parameter decreases from 50\% to 40\%. Note that,  as $\phi$ becomes small, the wind generation level becomes less volatile which results in a relatively low price volatility.   
This observation  indicates that the required storage capacity  to ensure a price volatility reduction target decreases as the wind power fluctuation parameter becomes small. Note that wind power fluctuation parameter $\phi$ can be reduced by improving the geographic diversity of wind farms in a region.

%The corresponding regulated and strategic storage capacities for achieving a maximum standard deviation of 1200 \$/MWh for wind power fluctuation of 40\% and 50\% are also compared in Fig. \ref{fig_vol}.

Based on Fig. \ref{fig_vol}, both price volatility and the required storage capacity for achieving a target price volatility become large as the wind power fluctuation parameter $\phi$ increases. 
%A strategic storage firm requires   60 and 70 MW more capacity, compared with a regulated storage firm, to achieve the  price volatility target (square root) of 1200 \$/MWh with wind power fluctuation parameter of 40\% and 50\%, respectively.
 To reduce the square root of price volatility to 1200 \$/MWh, the required strategic capacity with $\phi=40\%$ and $\phi=50\%$ is 60 MWh and 70 MWh, respectively, more than that of a regulated storage.
This observation confirms that regulated storage firms are more efficient  than strategic firms in reducing the price volatility. Although storage alleviates the price volatility in the market, it is not capable to eliminate it completely.

%Moreover, the optimum storage capacity required for achieving the price volatility target (square root) of 1200 \$/MWh indicates 60 and 70 MW more capacity requirement when the storage firms behave strategically for wind power fluctuation coefficients of 40 and 50\% respectively.   

\vspace*{-.2 cm}

\subsection{Two-node model simulations in South Australia and Victoria}

In the previous subsection, we analysed the impact of storage on the price volatility in SA when the SA-VIC interconnector is not active. In this subsection, we first study the  effect of  the interconnector between  SA and VIC on the price volatility in the absence of  storage firms. Next, we investigate the impact of storage firms on the price volatility when the SA-VIC transmission line operates at various capacities.
In our numerical results, SA is connected to VIC using a 680 MW interconnector which is currently operating with  70\% of its capacity, i.e. 30\% of its capacity is under maintenance.
The numerical results in this subsection are based on the two-node model for a high demand day with coal plant outage in SA. To investigate the impact of transmission line on  price volatility, it is assumed that the SA-VIC interconnector  operates with 60\% and 70\% of its capacity.

According to our numerical results, the peak price in SA and VIC is equal to 9634 \$/MWh when the SA-VIC interconnector is completely  in outage and  the wind power fluctuation parameter $\phi$ is equal to $50\%$. However, the peak price reduced to 1406 \$/MWh and 1114 \$/MWh when the interconnector operates at  60\% and 70\% of its capacity. The square root of price volatility is 2787 \$/MWh, 303 \$/MWh, and 219 \$/MWh when the capacity of the SA-VIC transmission line is equal to 0\%, 60\%, and 70\%, respectively. 

Simulation results show that as long as the transmission line is not congested, the interconnector alleviates the price volatility phenomenon in SA by importing electricity from VIC to SA at peak times. Since the market in SA compared to VIC is much smaller, about three times, the price volatility abatement in SA after importing electricity from VIC is much higher than the price volatility increment in VIC. Moreover, the price volatility reduces as the capacity of transmission line increases. 

\begin{figure}[!htb]  
	%\captionsetup{justification=centering}
	\centering
	\includegraphics[scale=0.6]{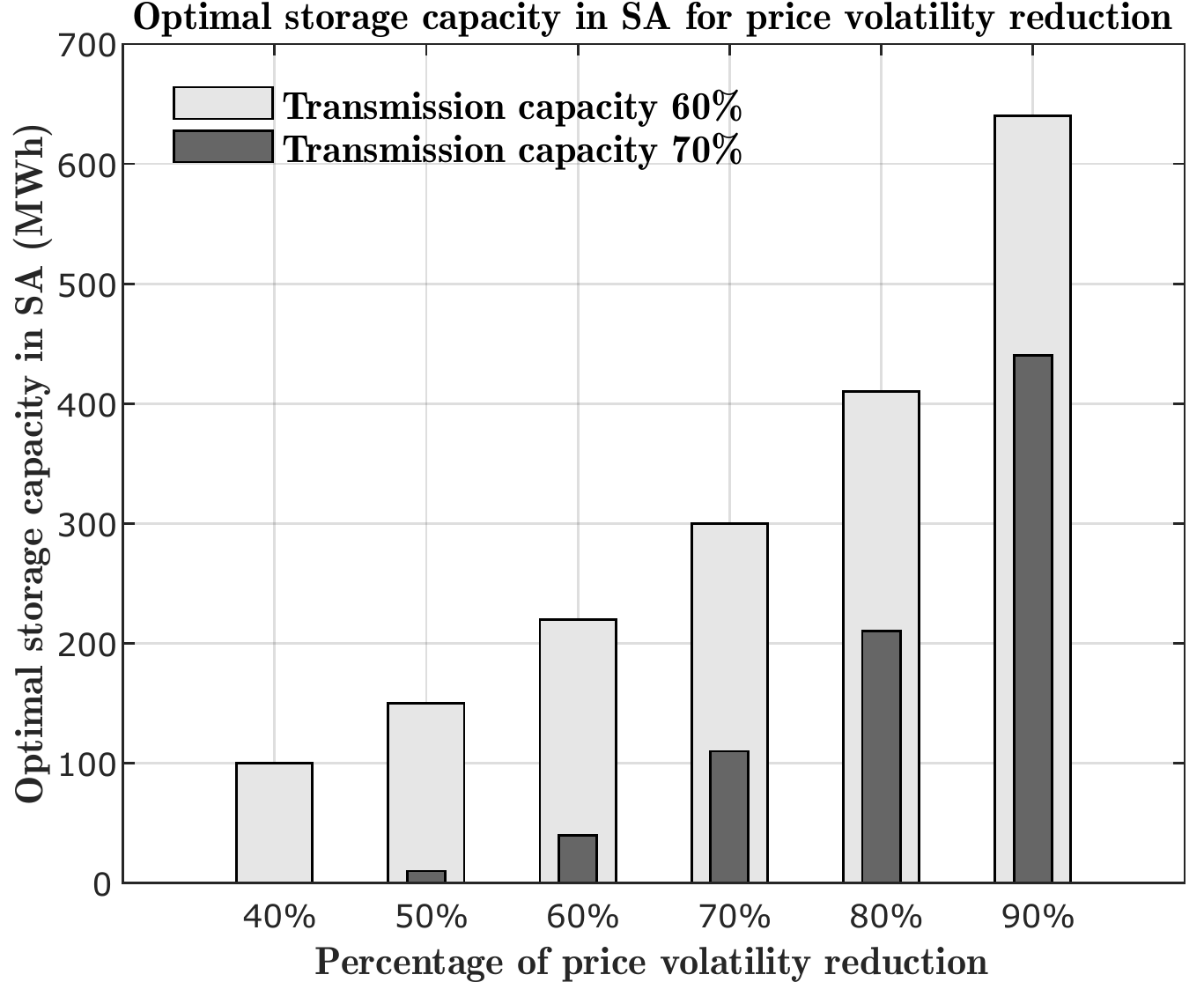} 
	\caption{Optimal regulated storage capacity versus the percentage of  price volatility reduction in the two-node market with wind power fluctuation parameter $\phi=50\%$ in a high demand day with coal-plants outage in SA.} 			
	\label{table2}
\end{figure}

%The wind power fluctuation with coefficient of 50\% creates peak price of 9634 \$/MWh when there is no interconnector. The 60\% and 70\%-capacity interconnectors reduce the peak price to 1406 and 1114 \$/MWh respectively. 
%Moreover, the higher the transmission line capacity, the lower the price volatility. The standard deviation of price is 2787, 303, and 219 \$/MWh for transmission capacity of 0\%, 60\% and 70\% respectively.
Fig. \ref{table2} shows the optimum storage capacity versus the percentage of price volatility reduction in the two-node market.  According to our numerical results,  storage is just located in SA, which witnesses a high level of  price volatility as the capacity of transmission line decreases. 
According to this figure, the optimum storage capacity becomes large as the capacity of transmission line decreases. Note that a sudden decrease of the transmission line capacity may result in a high level of price volatility in SA. However, based on Fig. \ref{table2}, storage firms are capable of reducing the price volatility during the outage of the   interconnecting lines.

\vspace*{-.2 cm}

\section{Conclusion} \label{sec: con}
High penetration of intermittent renewables, such as wind or solar farms, brings high levels of price volatility in electricity markets. Our study presents  an optimization model which decides on the minimum storage capacity required for achieving a price volatility  target in electricity markets. Based on our numerical results, the impact of storage on the  price volatility  in one-node electricity market of SA and two-node market of SA-VIC can be summarized as:

\begin{itemize}
\item Storage alleviates price volatility in the market due to the wind intermittency. However, storage does not remove price volatility completely, i.e. storage stops reducing the price volatility when it is not profitable.

\item The effect of a storage firm on  price volatility reduction depends on whether the firm  is  regulated or strategic. Both storage types have similar operation behaviour  and price reduction effects when they possess small capacities. For larger capacities, a strategic firm may under-utilize its available capacity and stop reducing the price level  due to its profit maximization strategy. On the other hand, a regulated storage firm is  more efficient  in price volatility reduction because of its social welfare maximization strategy. The price level and volatility reduction patterns observed when storage firms are regulated provide stronger incentives for the market operator to subsidize the storage technologies.

\item Both storage devices and transmission lines are capable of reducing the price volatility.  
High levels of price volatility that may happen due to the line maintenance  can be alleviated by storage devices.         

%\item Although we did not include the investment cost of storage technologies in our model, it is possible to compare the annual profit the storage firm gains from participating in the market and its technology cost to find out the rate of return for storage investment in the market. The $500 \rm{MWh}$ ($\frac{500}{2} \rm{MW}$) battery firm almost achieves a profit of 25 million\$/yr in our model and has the technology cost of 250 million\$ (with unit cost of 500 \$/kWh). Therefore, our roughly calculated results show that a battery storage firm returns its investment cost in almost 10 years.   

\end{itemize}

We intend to study the impact of ancillary services markets \cite{MulCommSpatial} and capacity markets \cite{CapEnergy}  on the integration of storage systems  in electricity networks in our future work. We also plan to investigate the impacts of other factors in both demand and supply sides on the price volatility.

\appendices

%\section{Inverse Demand Function} \label{App2}
%
%
%The linear and the exponential inverse demand functions crossing point (1500 MW, 100 \$/MWh) with the demand elasticity of price $\epsilon=10$ are visualized in Figure \ref{Exp_inverse}. We know that an unexpected load or storage charging shifts the inverse demand curve upward and results in higher price and demand with greater values for the demand elasticity. It can be seen that the linear inverse demand function decreases the elasticity ($\epsilon=2.9$) in such situations, \emph{i.e.} it does not capture the true relation between price and demand when storage facilities are charging. On the other hand, the exponential inverse demand function exhibits how an additional load or storage charging in the market leads to higher price, demand and elasticity ($\epsilon=10.6$).     
%\begin{figure}[!htb]  
%	%\captionsetup{justification=centering}
%	\centering
%	\includegraphics[scale= 0.65]{inverse.pdf} 
%	\caption{Comparison of linear and exponential inverse demand functions (base load \textthreequartersemdash, additional load \Kutline)} 			
%	\label{Exp_inverse}
%\end{figure} 

\vspace*{-.2 cm} 
\section{Charging/Discharging} \label{App1}
In this appendix, we show that the charge and discharge levels of any storage device  cannot be simultaneously positive at the NE of the lower game.
Consider a strategy in which  both charge and discharge levels of storage device $i$ at time $t$ under  scenario $w$, \emph{i.e.} $\qDis ,\qCh$, are both positive. We show that this strategy cannot be a NE strategy as follows. The net electricity flow of storage can be written as  $\qs=\etaDis \qDis - \frac{\qCh}{\etaCh}$. Let $\bar{q}^{\rm dis}_{itw}$ and $\bar{q}^{\rm ch}_{itw}$ be the new   discharge and charge levels of storage firm $i$ defined as $ \Big \{{{\bar{q}}^{\rm dis}_{itw}}=\qDis-\frac{\qCh}{\etaDis\etaCh}, \quad {{\bar{q}}^{\rm ch}_{itw}}=0 \Big \}$ if $ \qs>0 $, or $\Big \{{{\bar{q}}^{\rm dis}_{itw}}=0, \quad {{\bar{q}}^{\rm ch}_{itw}}=\qCh-\qDis \etaDis\etaCh \Big \}$ if $ \qs  <0 $.
%\begin{subequations}
%\begin{align*}
%&{{\bar{q}}^{\rm dis}_{itw}}=\qDis-\frac{\qCh}{\etaDis\etaCh}, \quad {{\bar{q}}^{\rm ch}_{itw}}=0 \quad {\rm if} \quad \qs>0\\
%&{{\bar{q}}^{\rm dis}_{itw}}=0, \quad {{\bar{q}}^{\rm ch}_{itw}}=\qCh-\qDis \etaDis\etaCh \quad {\rm if} \quad \qs  <0.
%\end{align*} 
%\end{subequations}
The new net flow of electricity can be written as  $\bar{q}^{\rm s}_{itw}=\etaDis {{\bar{q}}^{\rm dis}_{itw}}-\frac{{{\bar{q}}^{\rm ch}_{itw}}}{\etaCh}$. Note that the new variables $\bar{q}^{\rm s}_{itw}$, ${{\bar{q}}^{\rm ch}_{itw}}$ and ${{\bar{q}}^{\rm dis}_{itw}}$  satisfy the constraints
(\ref{cons41}-\ref{cons44}).

Considering the new charge and discharge strategies ${{\bar{q}}^{\rm dis}_{itw}}$ and ${{\bar{q}}^{\rm ch}_{itw}}$, instead of $\qDis$ and $\qCh$,  the nodal price and the net flow of storage device $i$ do not change. However, the charge/discharge cost of the storage firm $i$, under the new strategy, is reduces by:
	\begin{align}
 c_{i}^{\rm s} \left(\qCh+\qDis \right) - c_{i}^{\rm s} \left({\bar{q}}^{\rm dis}_{itw}+{\bar{q}}^{\rm ch}_{itw}\right) > 0  \nonumber
	\end{align}
Hence, any strategy in which the charge and discharge variables are simultaneously positive cannot be a NE, i.e. at the NE of the lower game each storage firm is either in the charge mode or discharge mode.

%\section{Monte Carlo simulation} \label{monCar}
%
%Within a Monte Carlo simulation, our model is simulated 500 times corresponding to random value of wind power availability ($\phi$). The effect of storage capacity (500 MWh) on electricity prices in SA in the High Peak Demand case is shown in Figure \ref{Exp_monte250}. It can be seen that the storage has significant impact on electricity prices which are very high.  
%
%
%\begin{figure}[!htb]  
%	%\captionsetup{justification=centering}
%	\centering
%	\includegraphics[scale= 0.65]{monte1.pdf} 
%	%\includegraphics[trim={3cm 9cm 1.5cm 9cm},clip,scale= 0.65]{monte1.pdf}
%	%\footnotesize{http://www.tnoconsulting.com.au/}
%	\caption{PDF of hourly electricity prices in SA in a high demand day with coal-plants outage case (500 samples of 24-hours in Monte Carlo simulation).} 			
%	\label{Exp_monte250}
%\end{figure}
%

\bibliographystyle{IEEEtran}
\bibliography{references}

%\pdftops

\end{document}